\documentclass[3p,times]{elsarticle}

\usepackage{ecrc}
\usepackage{amsmath,amssymb,amsthm}
\usepackage{hyperref}

\newtheorem{theorem}{Theorem}[section]

\newtheorem{proposition}{Proposition}[section]
\newtheorem{lemma}[theorem]{Lemma}
\newtheorem{definition}[theorem]{Definition}

\volume{00}

\firstpage{1}

\journalname{Differential Geometry and Its Applications}

\runauth{}

\jid{procs}

\CopyrightLine{2011}{Published by Elsevier Ltd.}

\usepackage{amssymb}
\usepackage[figuresright]{rotating}

\begin{document}

\begin{frontmatter}


 \title{Cartan meets Cramér-Rao}
 \author{Sunder Ram Krishnan}
 \ead{eeksunderram@gmail.com}
 %\fntext[label2]{}
 \address{Department of Computer Science and Engineering,
Amrita Vishwa Vidyapeetham, Amritapuri, 690525, Kerala, India}
 %\fntext[label3]{}

\dochead{}
%% e.g. \dochead{17th International Conference on Dynamical Processes in Excited States of Solids}

\title{Cartan meets Cramér-Rao}


\author{}

\address{}

\begin{abstract}
A Cartan--geometric, jet bundle formulation of curvature--aware variance bounds in parametric statistical estimation is developed. Building on our earlier extrinsic Hilbert space approach to the Cramér--Rao and Bhattacharyya--type inequalities, we show that the curvature corrections induced by the square root embedding of a statistical model admit a canonical intrinsic interpretation via jet geometry and Cartan’s prolongation theory.

For a scalar--parameter family with square root map
\( s_\theta=\sqrt{f(\cdot;\theta)}\in L^2(\mu) \),
we regard \( s_\theta \) as a section of the statistical bundle
\( E=\Theta\times L^2(\mu) \) and study its finite--order prolongations.
We point out that the classical algebraic efficiency condition—that the estimator residual
\( (T-\theta)s_\theta \) lies in the span of derivatives of \( s_\theta \) up to order \( m \)—is equivalent to the existence of a linear ordinary differential equation (ODE) of order \( m \) satisfied by the square root map.
Geometrically, this means the prolonged section lies in an ODE--defined submanifold of the jet bundle and is an integral curve of the restricted Cartan vector field.

The obstruction to such finite--order integrability is identified with the vertical component of the canonical Ehresmann connection on the jet tower, which coincides with the curvature correction term in variance bounds.
This establishes a direct correspondence between algebraic projection conditions in \( L^2(\mu) \) and intrinsic holonomy properties of statistical sections, yielding a unified geometric interpretation of higher--order information inequalities.
\end{abstract}


\end{frontmatter}

% \linenumbers

\section{Introduction}
\label{sec:intro}

The Cramér--Rao bound (CRB) and its higher--order analogs are central tools in
mathematical statistics for quantifying lower bounds on estimator variance, leading to the concept of estimator efficiency \cite{rao,vaart,schervish}. Notable classical generalizations of the CRB include the Chapman-Robbins bound \cite{chapman_robbins} and the Ziv-Zakai bound \cite{ziv_zakai}, which relax the usual regularity assumptions. Bayesian and hybrid variants of the CRB have also been developed \cite{hybrid_crb}. 
For CRB-inefficient estimators, another classical refinement is the \emph{Bhattacharyya inequality} (\cite{b1946}; see also \cite{ZS,l1998}) that uses higher-order derivatives of the log-likelihood. A Hilbert space foundation for the same can be found, for example, in \cite{pillai}; however, it lacks a geometric interpretation. 

Classical information geometry, pioneered by Amari et al., studies intrinsic affine
connections associated with the Fisher--Rao metric and the dual
$\alpha$--connections \cite{amari}; an important predecessor for this body of work was the influential paper of Efron \cite{efron1975}, who introduced ``statistical curvature" to quantify how close arbitrary one-parameter families are to exponential. It was also shown to be closely related to Fisher and Rao's theory of second-order efficiency. Turning to other works in information geometry, Smith \cite{smith_intrinsic_crb} derived intrinsic versions of the CRB on manifolds and a key follow-up was completed by Barrau and Bonnabel \cite{boumal2013_intrinsic_crb}. Extension of these results to the Bayesian setting may be found in Bouchard et al. \cite{bouchard2024} in the context of covariance matrix estimation, where performance bounds were derived for both the Euclidean metric and a natural affine invariant metric. 
Further, Li and Zhao \cite{Li} derived the Wasserstein--Cramér-Rao inequality, which is a lower bound on the Wasserstein variance of an estimator. Conditions for an estimator to be asymptotically efficient with respect to this bound were investigated by Nishimori and Matsuda \cite{nishimori2025}. 

There are also divergence-based generalizations of the CRB. A lower bound on the variance of an unbiased estimator for the $\alpha$-escort distribution was obtained by Mishra and Kumar in \cite{mishra2020}. A unifying differential geometric framework for deriving CRB-like inequalities from generalized Csiszár divergences was provided by the same authors in \cite{ashok2020} using Eguchi's theory. Further, relying on the Amari-Nagoaka framework, a Cramér-Rao type inequality was established in \cite{mishra2021}. All of these approaches typically remain in an intrinsic, divergence-induced geometry.

In recent work \cite{SRK,srkvector}, we developed an {\em extrinsic} geometric viewpoint
on the CRB by embedding the statistical model via the square root map
\[
s_\theta(\cdot)=\sqrt{f(\cdot;\theta)}\in L^2(\mu),
\]
viewing the parametric family 
\(\mathcal M=\{s_\theta:\theta\in\Theta\}\) as an immersed submanifold of the
ambient Hilbert space \(L^2(\mu)\).  
In that formulation, the first--order CRB corresponds to projecting the estimator
residual onto the tangent space spanned by the first derivative with respect to \(\theta\), $\operatorname{span}\{s_\theta^{(1)}\}$, while
classical higher--order bounds of the type of Bhattacharyya \cite{b1946}/Pillai--Sinha \cite{pillai} arise from projections onto the finite jet spans 
$\operatorname{span}\{s_\theta^{(1)},\dots,s_\theta^{(m)}\}$.  
Departures of the estimator error into directions normal to these jet spans
produce curvature corrections expressed through the second fundamental form of
the embedding $\mathcal M\subset L^2(\mu)$; quantifying this residual yields
sharpened variance lower bounds (see \cite{SRK,srkvector}).

The purpose of the present paper is to place this extrinsic programme on a
canonical {\em intrinsic} footing using the language of jet bundles and Cartan's
prolongation formalism -- the importance of the square root map in both the extrinsic and intrinsic pictures is also highlighted. For the former, this stems from the fact that we need a metric compatible connection whereas, for the latter, we need fixed fibre spaces. We show that the algebraic projection conditions underlying the CRB and higher--order Bhattacharyya-type inequalities
admit an equivalent description as differential equation (integrability/holonomy) conditions for the prolonged square root section in the statistical jet bundle. Recent work has explored Cartan–type metrics and invariant connections in information geometry~\cite{diatta}, yet their focus is completely different.

To reiterate and summarize, classical works in information geometry use constructions that do not rely on an ambient Hilbert space, contrasting with our earlier work \cite{SRK,srkvector} that introduced an extrinsic Hilbert-space embedding via the square root and derived curvature corrections to CRB-type bounds using the second fundamental form.  
Jet bundles and Cartan's prolongation formalism are classical in the geometry of
differential equations (e.g.\ \cite{Saunders}), but have not previously been
linked to variance bounds or estimator efficiency.  
The present paper fills this gap by providing a canonical intrinsic
(jet-theoretic) interpretation of the curvature terms obtained in
\cite{SRK}, thereby connecting classical CRB theory to Cartan geometry. We observe that the non--integrability of the Cartan distribution in finite jet order, along
with the canonical Ehresmann splitting of successive projections
in the jet bundle hierarchy, produces intrinsic geometric objects---specifically,
the vertical component of the canonical connection---that correspond exactly to
the curvature terms identified extrinsically in \cite{SRK}.  
This viewpoint not only supplies a coordinate--free explanation of the 
``second fundamental vector'' that appears for $m>1$ in \cite{SRK}, but also
organises the discussion on efficiency of unbiased estimators in the hierarchy of variance bounds as a sequence of holonomy/integrability constraints in the tower of higher jet bundles.

To facilitate further discussion and to avoid ambiguity, we state our assumptions and notation at the outset. 
\medskip

\noindent\textbf{Standing assumptions and notation.}
We adopt the following conventions throughout (consistent with \cite{SRK}):
\begin{itemize}
  \item The base parameter space is $\Theta\subset\mathbb R$.  
 The parametric family of distributions $\{P_\theta\}$ is absolutely continuous with
  respect to $\mu$, with densities $f(\cdot;\theta)\in \mathcal{H}$ and $\mathcal{H}=L^2(\mu)$, such that  
  $s_\theta=\sqrt{f(\cdot;\theta)}$ is $C^{m+1}$ in $\theta$ as an $\mathcal{H}$--valued
  map for each finite $m$ considered.

  \item For finite jet spans  
  $T_m(\theta)=\operatorname{span}\{s_\theta^{(1)},\dots,s_\theta^{(m)}\}$,  
  we assume the relevant Gram matrices are nonsingular, ensuring well-posed
  linear projection.
  
  \item Differentiation under the integral sign is permitted; with \(T(\cdot): \mathbb E_\theta[T]<\infty\),
\[
\frac{\partial}{\partial\theta}\int T(x) f(x;\theta)d\mu=\int T(x)\frac{\partial f(x;\theta)}{\partial\theta} d\mu,
\] 
whenever the right hand side is finite.
\end{itemize}
While the first assumption is important for us throughout the paper, the last two will be used in developing the statistical implications in Section \ref{sec:applications}.
\paragraph{Notations.}The statistical bundle is the trivial bundle $E=\Theta\times \mathcal{H}$. We write $J^m(E)$ for the $m$--jet bundle, $\omega_k$ for the canonical
  contact forms, and $D_\theta^{(m)}$ for the total derivative (Cartan) vector
  field on $J^m(E)$.  The $m$--jet prolongation of a section $s$ is denoted
  $j^m s$. Orthogonal projection onto a space \(W\subset \mathcal{H}\) is denoted by \(\Pi_W(\cdot)\).

With the setup behind us, we seek to clarify our contributions in more detail. 
\paragraph{Main contributions.}
\begin{enumerate}
  \item \emph{Cartan--jet reformulation of projection conditions.}  
  We point out that the algebraic membership condition
  \[
    (T-\theta)s_\theta \in \operatorname{span}\{s_\theta^{(1)},\dots,s_\theta^{(m)}\},
  \]
  which encodes first--order CRB efficiency ($m=1$) and higher--order
  Bhattacharyya/Pillai--Sinha-type efficiency, is equivalent to the existence of a
  linear differential relation of order $m$ satisfied by the square root map
  $s_\theta$.  
  Equivalently, the $m$-th prolongation $j^m s:\Theta\to J^m(E)$ takes values in 
  an ordinary differential equation (ODE)--defined submanifold $\Sigma_m\subset J^m(E)$ and is an integral curve
  of the restriction of the total derivative vector field.

  \item \emph{Intrinsic geometric source of curvature corrections.}  
  Working with the canonical contact forms and total derivative on $J^m(E)$ (for
  the statistical bundle $E=\Theta\times \mathcal{H}$), we identify the structural
  curvature/torsion obstructing integrability of the Cartan distribution on
  $J^{m}$.  
  Passing to $J^{m+1}$ removes this obstruction.  
  The corresponding vertical component of the Ehresmann connection on
  the affine bundle maps—via a canonical operator—into the
  $L^2$ residual of $(T-\theta)s_\theta$ orthogonal to $T_m(\theta)$, recovering
  exactly the curvature correction term present in the extrinsic Hilbert space
  picture.

  \item \emph{A dictionary between algebraic and geometric viewpoints.}  
  We establish a precise correspondence  \small
  \[
    \text{algebraic projection conditions in $L^2$}
    \Longleftrightarrow
    \text{holonomy/integrability conditions in the jet tower},
  \]\normalsize
  thereby explaining efficiency with respect to the CRB and higher--order Bhattacharyya-type bounds as
  successive holonomy constraints on the prolonged square root section.
  
  \item \emph{Conceptual clarification and further directions.}  
  We describe how the canonical jet torsion relates to the
  second fundamental form in the extrinsic picture, point out that the square root embedding provides a fixed fibre space so that the jet geometry applies
directly, and outline, among other directions, tensorial generalisations connecting jet geometry to intrinsic Fisher–Rao structures.
\end{enumerate}

\medskip
\noindent\textbf{Organization of the paper.}
Section~\ref{sec:prelim} reviews jet bundle preliminaries, Cartan distributions, structural curvature, the canonical Ehresmann connection,
horizontal/vertical decomposition, torsion, and a second fundamental form analogy. 
Section~\ref{sec:applications} establishes the main formal results with statistical implications: the  equivalence
between projection and ODE--integrability conditions, further linking the intrinsic jet bundle geometric development to the extrinsic curvature corrections of
\cite{SRK,srkvector}.  
Section~\ref{sec:examples} seeks to illustrate these geometric ideas with two concrete examples.  
We conclude in Section~\ref{sec:future} with a short discussion on future work.  

\medskip
The presentation is self-contained; readers familiar with jet geometry may skim through
Section~\ref{sec:prelim}.

% (References to \cite{KR1,KR2,Saunders,amari} match those in the bibliography.)

\section{Jet bundle preliminaries and Cartan prolongations}
\label{sec:prelim}

In this section we recall the basic geometry of jet bundles needed in the rest of the paper.  
Our presentation follows the standard treatment in \cite{Saunders}, but is specialised to the scalar base
variable and the trivial bundle structure relevant for statistical applications.  

\subsection{Jets of sections for a scalar base variable}

Let $E = \mathbb{R} \times \mathcal{H}$ denote the trivial fibre bundle over the base parameter manifold 
$\Theta = \mathbb{R}$ with projection $\pi:E\to \Theta$, $\pi(\theta,s)=\theta$.  
A smooth section $s:\Theta\to E$ assigns to each $\theta$ a fibre value $s(\theta)$ and will be denoted
simply by $s(\theta)$ or $s_\theta$.  

The $m$-jet bundle $J^{m}(E)$ encodes equivalence classes of sections that agree up to $m$-th order in derivatives.
In coordinates a point of $J^{m}(E)$ is represented by
\[
 (\theta, s, s_1, s_2,\ldots,s_m),
\]
where $s_k$ records the value of 
$s^{(k)}(\theta)$.  
The canonical projections
\[
\pi_{m,0}:J^{m}(E)\to E,\qquad
\pi_{m,m-1}:J^{m}(E)\to J^{m-1}(E)
\]
forget higher-order coordinates.

The $m$-jet prolongation of a section $s$ is the map
\[
 j^{m}s(\theta)=\big(\theta, s(\theta), s^{(1)}(\theta),\ldots,s^{(m)}(\theta)\big).
\]
A smooth curve in $J^{m}(E)$ is \emph{holonomic} if it arises as $j^{m}s$ for some $s$; holonomicity is encoded
by the vanishing of the contact forms below.

\subsection{Cartan distribution and total derivative}

The bundle $J^{m}(E)$ carries canonical contact $1$-forms
\[
 \omega_k = ds_k - s_{k+1}\,d\theta,\qquad k=0,\ldots,m-1.
\]
These vanish on all prolonged curves $j^m s$ since $ds_k = s^{(k+1)} d\theta$ along such curves.

The {\em Cartan distribution} $\mathcal{C}^{m}\subset TJ^{m}(E)$ is the common kernel
\[
 \mathcal{C}^{m}=\bigcap_{k=0}^{m-1}\ker(\omega_k),
\]
and is $1$-dimensional.  Its generator is the {\em total derivative} (Cartan) vector field
\begin{equation}
\label{eq:total-derivative}
D^{(m)}_\theta
 = \frac{\partial}{\partial \theta} 
   + s_1\frac{\partial}{\partial s}
   + s_2\frac{\partial}{\partial s_1}
   + \cdots
   + s_m\frac{\partial}{\partial s_{m-1}},
\end{equation}
which is the restriction of the formal operator
\[
 D_\theta 
   = \frac{\partial}{\partial \theta}
     + \sum_{k\ge 0} s_{k+1}\frac{\partial}{\partial s_k}
\]
to $J^{m}(E)$.  
The field $D^{(m)}_\theta$ is the canonical horizontal lift of $\partial_\theta$ and represents the flow
which preserves holonomicity to order $m$.

The Cartan distribution has rank one, while the vertical bundle
\[
 VJ^{m}= \ker((\pi_{m,0})_*)
\]
has rank $m+1$.  Thus $TJ^{m}=\mathcal{C}^{m}\oplus VJ^{m}$ as a vector bundle.

For the projection $\pi_{m+1,m}:J^{m+1}\to J^{m}$, the relative vertical space is 
\[
 V(J^{m+1}/J^{m})
 = \ker((\pi_{m+1,m})_*)
 = \mathrm{span}\Big\{\frac{\partial}{\partial s_{m+1}}\Big\},
\]
corresponding to the new derivative coordinate $s_{m+1}$ introduced at level $m+1$.

\subsection{Geometry of prolonged curves}

Let $j^{m}s(\theta)$ be a prolonged curve.  
The following result makes precise the fundamental decomposition of the $\theta$-derivative of a prolonged section.

\begin{lemma}[Jet prolongation identity]
\label{lem:jet-prolong}
Let $s:\mathbb{R}\to E$ be a smooth section with $m+1$ derivatives.  
Then the $\theta$-derivative of the prolonged curve $j^{m}s$ is given by
\begin{equation}
\label{eq:jet-prolong-identity}
\frac{d}{d\theta}j^{m}s(\theta)
=
D^{(m)}_\theta
\;+\;
s^{(m+1)}(\theta)\,\frac{\partial}{\partial s_m},
\end{equation}
where $D^{(m)}_\theta$ is the Cartan total derivative vector field
on $J^{m}(E)$.
\end{lemma}

\begin{proof}
This is an immediate consequence of the coordinate expression of $j^{m}s$ and of the
definition of the Cartan vector field~\eqref{eq:total-derivative}.  
Differentiating the coordinate representation
\(
j^{m}s(\theta)=\big(\theta,s(\theta),s^{(1)}(\theta),\dots,s^{(m)}(\theta)\big)
\)
yields the stated decomposition.
\end{proof}

Thus the velocity of $j^{m}s$ decomposes into:
(i) a horizontal part $D^{(m)}_{\theta}$ spanning $\mathcal{C}^m$, and  
(ii) a vertical part $s^{(m+1)}(\theta)\,\frac{\partial}{\partial s_m}$.

The vertical component is precisely the correction required to ensure that the lifted curve becomes holonomic
one level higher, in $J^{m+1}(E)$, in order to satisfy the new contact form $\omega_m=ds_m - s_{m+1}d\theta$.

\subsection{Structural curvature and non-integrability}

The Cartan distribution $\mathcal{C}^{m}$ is not integrable for finite $m$.  
The obstruction is encoded in the exterior derivative of the highest-order contact form
\[
 \omega_{m-1} = ds_{m-1} - s_m d\theta.
\]
\begin{theorem}[Non-integrability of the Cartan distribution]
\label{thm:non-integrable}
Let $\mathcal{C}^{m}=\bigcap_{k=0}^{m-1}\ker(\omega_k)$ be the Cartan distribution on $J^{m}(E)$.
Then $\mathcal{C}^{m}$ is never involutive for any finite $m$.  
More precisely,
\[
d\omega_{m-1}=d\theta\wedge ds_{m}
\not\equiv 0 
\quad \text{mod }\{\omega_0,\dots,\omega_{m-1}\},
\]
so that $\mathcal{C}^{m}$ fails the Frobenius involutivity condition.
\end{theorem}

\begin{proof}
Since
$\omega_{m-1}=ds_{m-1}-s_m d\theta$, we have $d\omega_{m-1}=d\theta\wedge ds_m$, which
cannot be expressed in the ideal generated by $\{\omega_0,\dots,\omega_{m-1}\}$ on $J^{m}(E)$.
Hence the Cartan distribution is non-integrable by the Frobenius theorem.
\end{proof}

However, passing to $J^{m+1}(E)$ introduces the new contact form $\omega_m=ds_m-s_{m+1}d\theta$, so that
\[
 d\omega_{m-1}
   = d\theta\wedge \omega_m
   \equiv 0
   \qquad \text{mod } \{\omega_0,\ldots,\omega_m\}.
\]
Thus the obstruction to integrability at order $m$ becomes absorbed into the contact ideal at order $m+1$.
This phenomenon recurs throughout the jet tower and underlies the prolongation method for differential
equations.

\subsection{Ehresmann connection and decomposition}

The projection $\pi_{m,m-1}:J^{m}(E)\to J^{m-1}(E)$ is an affine bundle.
The {\em canonical Ehresmann connection} $\nabla^{m}$ is defined by choosing the horizontal subbundle
\[
 H^{m}
   = \mathrm{span}\big\{
      D^{(m)}_{\theta},
      \tfrac{\partial}{\partial s},
      \tfrac{\partial}{\partial s_1},
      \ldots,
      \tfrac{\partial}{\partial s_{m-1}}
    \big\}.
\]
Since $V^{m}\equiv V(J^{m}/J^{m-1})=\mathrm{span}\{\partial/\partial s_m\}$,
we obtain a direct-sum splitting
\[
 TJ^{m}=H^{m}\oplus V^{m}.
\]

For the tangent vector
\[
 X^{(m)}:=\frac{d}{d\theta}j^{m}s
 = D^{(m)}_{\theta}
   + s^{(m+1)}\frac{\partial}{\partial s_m},
\]
the horizontal and vertical projections are
\[
 X^{(m)}_H = D^{(m)}_{\theta},\qquad
 X^{(m)}_V = s^{(m+1)}\frac{\partial}{\partial s_m},
\]
respectively.
Lemma \ref{lem:jet-prolong} already captured the geometric content: the component along $\partial/\partial s_m$ is
\emph{forced} by the requirement that the curve remain holonomic upon passing from $J^{m}$ to $J^{m+1}$. With the connection perspective, we have:

\begin{proposition}[Horizontal/vertical decomposition]
\label{prop:horiz-vert}
For the canonical Ehresmann connection $\nabla^{m}$ on the affine bundle 
$\pi_{m,m-1}:J^{m}(E)\to J^{m-1}(E)$, the tangent vector in~\eqref{eq:jet-prolong-identity} decomposes as
\[
\frac{d}{d\theta} j^{m}s
=
\underbrace{D^{(m)}_\theta}_{X^{(m)}_H}
\;+\;
\underbrace{s^{(m+1)}(\theta)\frac{\partial}{\partial s_m}}_{X^{(m)}_V},
\qquad
X^{(m)}_H\in H^{m},\ \ X^{(m)}_V\in V^{m}.
\]
This decomposition is unique.
\end{proposition}

\begin{proof}
By construction $H^{m}$ is spanned by $D^{(m)}_\theta$ and $\partial/\partial s_k$ for $k< m$,
whereas $V^{m}=\ker((\pi_{m,m-1})_*)=\mathrm{span}\{\partial/\partial s_m\}$ is one-dimensional.
The decomposition is immediate.
\end{proof}

The vertical component $X^{(m)}_V$ measures the obstruction to remaining within the horizontal distribution
$H^{m}$.  
Equivalently, $X^{(m)}_V$ is the obstruction to flatness: it is the ``acceleration’’ required to satisfy the
new contact form on $J^{m+1}(E)$.

The {\em torsion form} of the canonical connection is obtained by contracting the structural curvature 
$d\omega_{m-1}$ with the base total derivative:
\[
 \Omega_m = i_{D_\theta} d\omega_{m-1}.
\]
In Theorem \ref{thm:non-integrable}, we saw that when passing to the next jet bundle level, integrability is restored modulo the higher contact forms.  
This observation is encoded by the torsion form of the canonical connection.

\begin{proposition}[Torsion form]
\label{prop:torsion}
For the canonical connection $\nabla^{m}$, the torsion form is
\[
\Omega_m=i_{D_\theta}d\omega_{m-1}
     =\omega_{m},
\]
the highest-order contact form on $J^{m+1}(E)$.
\end{proposition}

\begin{proof}
Applying the interior product $i_Y (\alpha \wedge \beta) = (i_Y \alpha) \beta - \alpha (i_Y \beta)$:
\[
\Omega_m = d\theta(D_\theta) ds_m - ds_m(D_\theta) d\theta=ds_m-s_{m+1}d\theta. 
\]But $ds_m-s_{m+1}d\theta$ is $\omega_m$ on $J^{m+1}(E)$.
\end{proof}

Thus the vertical component in Proposition~\ref{prop:horiz-vert} is
the necessary geometric correction that cancels torsion when passing from $J^{m}$ to $J^{m+1}$, reflecting again that the vertical component $X_V$ encodes the failure of $j^{m}s$ to be horizontal in $J^{m}$.
This is the jet-analog of the normal component of acceleration in the second fundamental form of a
curve.

\subsection{Second fundamental form analogy}

Interpreting $J^{m}(E)\subset J^{m+1}(E)$ as an inclusion of manifolds,
the vertical component $X^{(m)}_V$ of the jet velocity plays the role of a ``normal acceleration,’’ analogous to
the second fundamental form of a curve.  
The horizontal component $X^{(m)}_H$ provides the ``tangential acceleration’’ determined by the Cartan flow
$D^{(m)}_\theta$.
This analogy anticipates the comparison between intrinsic jet geometry and the extrinsic Hilbert space curvature
used in \cite{SRK}.  
We return to this interpretation in Section~\ref{sec:applications} when discussing statistical
implications.

%%%%%%%%%%%%%%%%%%%%%%%%%%%%%%%%%%%%%%%%%%%%%%%%%%%%%%%%%%%%%%%%%%%%%%%%%%%%%%
\section{Statistical Applications: Jet Geometry and Variance Bounds}
\label{sec:applications}
%%%%%%%%%%%%%%%%%%%%%%%%%%%%%%%%%%%%%%%%%%%%%%%%%%%%%%%%%%%%%%%%%%%%%%%%%%%%%%

We now translate the geometric identities of the preceding section into
statistical statements for variance bounds.  The key principle is:

\[
\boxed{
\text{Algebraic projection conditions in }L^2
\;\Longleftrightarrow\;
\text{ODE/holonomy constraints in jet space}.
}
\]
\subsection{Projection conditions and ODE--integrability}

Let
\[
T_m(\theta)=
\mathrm{span}\{s_\theta^{(1)},\dots,s_\theta^{(m)}\}\subset \mathcal{H}.
\]
For an unbiased estimator $T(X)$, the classical algebraic
projection condition characterising $m$-th order efficiency is
\begin{equation}\label{eq:proj-cond}
(T(X)-\theta)s_\theta \in T_m(\theta).
\end{equation}

Define the $m$-th order ODE--constraint submanifold
\begin{equation}
\label{eq:sm}
\Sigma_m
=
\bigl\{
(\theta,s,s_1,\dots,s_m)\in J^m(E):
(T-\theta)s=c_{1,\theta}s_1+\cdots+c_{m,\theta}s_m
\bigr\},
\end{equation}
where the coefficients $c_{k,\theta}$ are the coordinate coefficients of the
orthogonal projection of $(T-\theta)s_\theta$ onto $T_m(\theta)$ and we assume that the vectors in \(T_m\) are linearly independent. We have the following result, where the \(m=1\) case corresponds to CRB efficiency.

\begin{theorem}[Cartan viewpoint of $m$-th order efficiency]
\label{thm:equivalence}
For any $m\ge1$, the following are equivalent:
\begin{enumerate}
\item The projection condition \eqref{eq:proj-cond} holds for all $\theta$.
\item The prolonged section $j^m s(\theta)$ takes values in $\Sigma_m$ and
the total derivative $D_\theta^{(m)}$ restricted to $\Sigma_m$ is tangent to
$\Sigma_m$; i.e.
\[
D_\theta^{(m)}\big|_{j^ms(\theta)}\in T_{j^m s(\theta)}\Sigma_m.
\]
\end{enumerate}
\end{theorem}

\begin{proof}
In coordinates on $J^m(E)$,
\[
D_\theta^{(m)}=
\partial_\theta +
s_1\partial_s +
\cdots +
s_m\partial_{s_{m-1}}.
\]
The constraint defining $\Sigma_m$ is linear in $(s,s_1,\dots,s_m)$.  A
vector is tangent to $\Sigma_m$ iff it annihilates the constraint
one--form
\[
\alpha_m
=
(T-\theta)sd\theta - c_{1,\theta}ds - \cdots - c_{m,\theta}ds_{m-1}.
\]
Evaluating $\alpha_m$ on $D_\theta^{(m)}$ gives exactly the defect of
$(T-\theta)s_\theta$ from the span of $s_\theta^{(1)},\dots,s_\theta^{(m)}$.
Thus $\alpha_m(D_\theta^{(m)})=0$ iff the projection condition holds.
\end{proof}
Thus, the CRB is attained exactly when the estimator error is
parallel to the score function.  
In jet language, this means that the prolongation satisfies a first-order ODE constraint,
so that $j^{1}s_\theta$ is an integral curve in the corresponding submanifold of $J^{1}(E)$.
When the CRB is \emph{not} attained, the prolongation $j^{1}s$ fails to lie in any such ODE 
submanifold.  
This is exactly the geometric obstruction to remaining horizontal in $J^{1}(E)$.

Next, we endeavour to connect the intrinsic geometry, developed here, with the extrinsic curvature-aware corrections derived in \cite{SRK, srkvector} using the statistical residual defined as the orthogonal projection of the estimator error onto \(T_j^\perp\):
\[
R^{(j)}_\theta
:=
\Pi_{T_j^\perp}((T-\theta)s_\theta)
\in \mathcal{H}.
\]
%%%%%%%%%%%%%%%%%%%%%%%%%%%%%%%%%%%%%%%%%%%%%%%%%%%%%%%%%%%%%%%%%%%%%%%%%%%%%%
\subsection{Jet decomposition and the statistical residual}
%%%%%%%%%%%%%%%%%%%%%%%%%%%%%%%%%%%%%%%%%%%%%%%%%%%%%%%%%%%%%%%%%%%%%%%%%%%%%%
Consider an $m$-th order efficient estimator (\(m> 1\)). From Theorem \ref{thm:equivalence}, we know that \((T-\theta)s=c_{1,\theta}s_1+\cdots+c_{m,\theta}s_m\), with \(c_{m,\theta}\neq 0\), for a unique set of coefficients \(\{c_{k,\theta}\}\). Also, from Proposition \ref{prop:horiz-vert} and the theorem again, we see that the prolonged derivative decomposes in jet coordinates as
\begin{eqnarray*}
X^{(m-1)}&:=&\frac{d}{d\theta}j^{m-1} s\\
&=&
D_\theta^{(m-1)} + X_V^{(m-1)},
\quad
X_V^{(m-1)} = \frac{(T-\theta)s-c_{1,\theta}s_1-\cdots+c_{m-1,\theta}s_{m-1}}{c_{m,\theta}}\partial_{s_{m-1}},
\end{eqnarray*}
when restricted to \(\Sigma_m\).
The key point, as we note below, is that
\(X_V^{(m-1)}\) captures the jet-theoretic obstruction to 
\((m-1)\)-th order efficiency.

To link the intrinsic and extrinsic geometries alluded to before, we define an operator that transfers vertical jet information into Hilbert-space residual information.

\begin{definition}[Jet verticals and residual]
\label{prop:transfer}
For each $\theta$ and fixed $m\geq 2$, define a map
\[
A_\theta:
V_{j^{m-1}s(\theta)}(J^{m-1}/J^{m-2})
\longrightarrow 
T_{m-1}(\theta)^\perp\subset \mathcal{H}
\]
such that
\[
R^{(m-1)}_\theta = A_\theta\bigl(X_V^{(m-1)}(\theta)\bigr),
\]
by
\[
A_\theta(s_m\partial_{s_{m-1}}) 
:= 
\Pi_{T_{m-1}^\perp}(s_\theta^{(m)}),
\]
\end{definition}

%%%%%%%%%%%%%%%%%%%%%%%%%%%%%%%%%%%%%%%%%%%%%%%%%%%%%%%%%%%%%%%%%%%%%%%%%%%%%%
%\subsection{Intrinsic Bhattacharyya identity}
%%%%%%%%%%%%%%%%%%%%%%%%%%%%%%%%%%%%%%%%%%%%%%%%%%%%%%%%%%%%%%%%%%%%%%%%%%%%%%
What we argued in \cite{SRK,srkvector} may then be stated as follows: Given an \((m-1)\)-th order inefficient estimator, we can derive curvature-aware corrections to CRB/Bhattacharyya-type bounds. Specifically,
\begin{equation}
\mathrm{Var}_\theta(T)
\geq 
\left\|\Pi_{T_{m-1}}\bigl((T-\theta)s_\theta\bigr)\right\|_{\mathcal{H}}^2
+
\left\|\Pi_{R^{(m-1)}_\theta}\bigl((T-\theta)s_\theta\bigr)\right\|_{\mathcal{H}}^2.
\label{eq:var-decomp}
\end{equation}
Note that:
\begin{enumerate}
\item The first term in \eqref{eq:var-decomp} is an $(m-1)$th-order Bhattacharyya-type bound.
\item The second term is nonnegative and vanishes iff $T$ is $(m-1)$-th order efficient. Thus we obtain a strict improvement whenever we have an inefficient estimator.
\item The vertical jet component $X_V^{(m-1)}$, mapped via \(A_\theta\), captures the correction to classical variance bounds along the direction of the second fundamental form viewed from the intrinsic perspective of the statistical jet bundle geometry \(J^{m-1}(E)\subset J^m(E)\). This was (slightly informally) introduced as a ``second fundamental vector" in our extrinsic development \cite{SRK}.
\end{enumerate}

\subsection{Importance of the square root embedding}
It is pertinent to stress the importance of the square root map to both the extrinsic and intrinsic pictures we developed. Just as we noted in \cite{SRK,srkvector} that the fixed Hilbert space \(\mathcal{H}\) provides access to a metric-compatible Levi-Civita connection in the extrinsic geometric perspective, having a fixed fibre space here is central to the intrinsic geometric development. In the jet bundle geometry, we need the fibres to be independent of the base parameter \(\theta\) precluding the space \(L^2(P_\theta)\). Therefore, working with the log-likelihood embedding as in the usual Bhattacharyya-type bounds, for instance, will not be feasible if one desires a clean geometric interpretation. In other words, we do not get a systematic (intrinsic/extrinsic) geometric explanation of the additional, non-asymptotic terms in classical variance bounds as curvature-aware corrections in the event of inefficiency.

\section{Examples and explicit computations}
\label{sec:examples}

We present two detailed examples that illustrate the results discussed.
The first is the Gaussian location model, where all jet quantities and inner products are computable in closed form;
the second treats a one-parameter exponential family to demonstrate a nontrivial torsion/curvature correction.

\subsection{Gaussian location model}
\label{sec:gauss}

Consider the scalar Gaussian density with known variance $\sigma^2>0$:
\[
f(x;\theta)=\frac{1}{\sqrt{2\pi}\sigma}\exp\!\Big(-\frac{(x-\theta)^2}{2\sigma^2}\Big),\qquad x\in\mathbb{R},\ \theta\in\mathbb{R}.
\]
The square root embedding is
\[
s_\theta(x)=\Big(\frac{1}{\sqrt{2\pi}\sigma}\Big)^{1/2}
            \exp\!\Big(-\frac{(x-\theta)^2}{4\sigma^2}\Big).
\]
Below we compute derivatives, inner products and show the ODE/projection condition for the first order.

Differentiating under the integral sign (justified by the Gaussian tails), one obtains
\[
s_\theta^{(1)}(x)=\frac{\partial}{\partial\theta}s_\theta(x)
                 =\frac{x-\theta}{2\sigma^2}\,s_\theta(x).
\]

All inner products are computed with respect to the Lebesgue measure (equivalently expectations under $P_\theta$ using $s_\theta^2=f(\cdot;\theta)$). Using Gaussian moments about the mean $\theta$ we get
\begin{align*}
\langle s_\theta, s_\theta\rangle &= \int s_\theta^2(x)\,dx = 1,\\
\langle s_\theta, s_\theta^{(1)}\rangle &= \tfrac12\frac{d}{d\theta}\int s_\theta^2(x)\,dx = 0,\\
\langle s_\theta^{(1)}, s_\theta^{(1)}\rangle
&= \int \frac{(x-\theta)^2}{4\sigma^4}s_\theta^2(x)\,dx
= \frac{1}{4\sigma^4}\,\operatorname{Var}_\theta(X)
= \frac{1}{4\sigma^2}.
\end{align*}
Hence the Fisher information (in the square root embedding convention \(I=4\langle s^{(1)},s^{(1)}\rangle\)) is
\[
I(\theta)=4\langle s_\theta^{(1)},s_\theta^{(1)}\rangle=\frac{1}{\sigma^2},
\]
which is the well--known result for the Gaussian location model.

Consider the unbiased estimator $T(X)=X$. The estimator error times $s_\theta$ is
\[
(T-\theta)s_\theta = (X-\theta)s_\theta.
\]
Compare with the score direction:
\[
s_\theta^{(1)}(x)=\frac{x-\theta}{2\sigma^2}s_\theta(x).
\]
Thus pointwise (for each $x$)
\[
(T-\theta)s_\theta = 2\sigma^2\, s_\theta^{(1)}.
\]
This verifies the first--order ODE/projection condition of Theorem~\ref{thm:equivalence} and satisfies \eqref{eq:sm} with $m=1$, $c_{1,\theta}=2\sigma^2$, and hence $T$ attains the CRB (indeed $X$ is efficient in this model). In our jet language, $j^1 s_\theta$ lies in the first--order ODE submanifold.

\subsection{Natural exponential family}
\label{sec:expfam-example}

Let
\[
f(x;\theta)=\exp\{\theta x - A(\theta)\}\,h(x),\qquad \theta\in\Theta,
\]
be a regular one--parameter natural exponential family with $A$ smooth and
$A''(\theta)>0$ for all $\theta\in\Theta$.  
Denote central moments by
\[
\mu_r(\theta):=\mathbb{E}_\theta\big[(X-A'(\theta))^r\big],\qquad r\ge1,
\]
so that in particular $\mu_2(\theta)=A''(\theta)$.

\paragraph{Square root derivatives}
The square root map and its first two derivatives are
\[
s_\theta(x)=\exp\!\Big(\tfrac12(\theta x - A(\theta))\Big)\sqrt{h(x)},
\]
\[
s_\theta^{(1)}(x)=\tfrac12\,(x-A'(\theta))\,s_\theta(x),
\qquad
s_\theta^{(2)}(x)=\Big(\tfrac14(x-A'(\theta))^2-\tfrac12 A''(\theta)\Big)s_\theta(x).
\]

\paragraph{Inner products and Gram matrix}
Using the identification
\[
\langle \Phi s_\theta,\Psi s_\theta\rangle
=\mathbb{E}_\theta[\Phi(X)\Psi(X)],
\]
we compute
\begin{align*}
\langle s_\theta^{(1)},s_\theta^{(1)}\rangle
&=\tfrac14\,\mathbb{E}_\theta[(X-A')^2]
=\tfrac14\,\mu_2,\\[4pt]
\langle s_\theta^{(1)},s_\theta^{(2)}\rangle
&=\mathbb{E}_\theta\!\Big[\tfrac12(X-A')\Big(\tfrac14(X-A')^2-\tfrac12A''\Big)\Big]
=\tfrac18\,\mu_3,\\[4pt]
\langle s_\theta^{(2)},s_\theta^{(2)}\rangle
&=\mathbb{E}_\theta\!\Big[\Big(\tfrac14(X-A')^2-\tfrac12A''\Big)^2\Big]
=\tfrac1{16}\,\mu_4,
\end{align*}
where, in the last identity, the cross terms involving $\mu_2=A''$ cancel exactly.

Hence the $2\times2$ Gram matrix for the jet basis
$\{s_\theta^{(1)},s_\theta^{(2)}\}$ is
\[
G_\theta
=
\begin{pmatrix}
\frac14\mu_2 & \frac18\mu_3\\[4pt]
\frac18\mu_3 & \frac1{16}\mu_4
\end{pmatrix},
\qquad
\det G_\theta=\frac1{64}\big(\mu_2\mu_4-\mu_3^2\big).
\]

\paragraph{Orthogonal projection onto the second--order jet span}
Fix an estimator $T=T(X)$.  
We decompose
\begin{equation}
(T-\theta)s_\theta
=
a(\theta)\,s_\theta^{(1)}
+
b(\theta)\,s_\theta^{(2)}
+
\mathcal{R}_\theta,
\qquad
\mathcal{R}_\theta\perp\operatorname{span}\{s_\theta^{(1)},s_\theta^{(2)}\},
\label{eq:exp-proj}
\end{equation}
where $a(\theta),b(\theta)$ are the coefficients of the
\emph{$L^2$--orthogonal projection} of $(T-\theta)s_\theta$ onto the
second--order jet span and $\mathcal{R}_\theta$ is the residual.

Taking inner products of \eqref{eq:exp-proj} with
$s_\theta^{(1)}$ and $s_\theta^{(2)}$ yields the linear system
\[
G_\theta
\begin{pmatrix} a(\theta) \\[2pt] b(\theta)\end{pmatrix}
=
\begin{pmatrix}
r_1(\theta)\\[2pt]
r_2(\theta)
\end{pmatrix},
\]
where
\[
r_1(\theta)
:=\langle (T-\theta)s_\theta, s_\theta^{(1)}\rangle
=\tfrac12\,\mathbb{E}_\theta[(T-\theta)(X-A')],
\]
\[
r_2(\theta)
:=\langle (T-\theta)s_\theta, s_\theta^{(2)}\rangle
=\mathbb{E}_\theta\,\!\Big[(T-\theta)\Big(\tfrac14(X-A')^2-\tfrac12A''\Big)\Big].
\]

Assuming $\det G_\theta\neq0$, the projection coefficients are uniquely determined as
\begin{align}
a(\theta)
&=
\frac{4\,\mu_4\,r_1 - 8\,\mu_3\,r_2}{\mu_2\mu_4-\mu_3^2},
\nonumber\\[4pt]
b(\theta)
&=
\frac{-8\,\mu_3\,r_1 + 16\,\mu_2\,r_2}{\mu_2\mu_4-\mu_3^2}.
\label{eq:ab-sol}
\end{align}

\paragraph{Second--order efficiency}
It is easily checked that the first--order efficiency condition does not hold here in general; more on this in our comments below. Second--order efficiency requires the stronger condition
\[
\mathcal{R}_\theta\equiv0,
\]
i.e.\ that $(T-\theta)s_\theta$ lie \emph{exactly} in the finite jet span given by 
$\operatorname{span}\{s_\theta^{(1)},s_\theta^{(2)}\}$.
While the $L^2$--projection \eqref{eq:exp-proj} always exists (when
$\det G_\theta\neq0$), the vanishing of $\mathcal{R}_\theta$ imposes strong
distributional constraints and does not generically hold in exponential families.
Suppose there exists an unbiased estimator $T$ such that
\[
(T-\theta)s_\theta
=
a(\theta)s_\theta^{(1)} + b(\theta)s_\theta^{(2)}
\quad\text{pointwise a.e. in \(x\)}
\]
for all $\theta\in\Theta$.
Dividing by $s_\theta(x)>0$, this is equivalent to the identity
\begin{equation}
T(x)-\theta
=
\alpha(\theta)\,(x-A'(\theta))
+
\beta(\theta)\big((x-A'(\theta))^2-2A''(\theta)\big),
\label{eq:quad-identity}
\end{equation}
where $\alpha(\theta)=\tfrac12 a(\theta)$ and
$\beta(\theta)=\tfrac14 b(\theta)$.

Thus second--order efficiency forces the unbiased estimator $T(X)$ to satisfy moment constraints
that constitute a strong structural restriction.
Except in special cases, such as Gaussian location families with usual estimators, where
$s_\theta^{(1)}$ already spans $(T-\theta)s_\theta$ pointwise,
the quadratic identity \eqref{eq:quad-identity} fails and the residual
$\mathcal{R}_\theta$ in \eqref{eq:exp-proj} is nonzero.

\paragraph{Remarks and special cases}
\begin{enumerate}
\item If $T$ is unbiased for \(\theta\) (i.e.\ \(\mathbb{E}_\theta[T]=\theta\)),
then $\mathbb{E}_\theta[T-\theta]=0$, the formulae simplify, and we may write
\[
r_1(\theta)=\tfrac12\,\mathrm{Cov}_\theta(T,X),\qquad
r_2(\theta)=\tfrac14\,\mathbb{E}_\theta[(T-\theta)(X-A')^2].
\]
Note that these depend on $\theta$ and on the chosen estimator $T$.
\item For a fixed $\theta$ and a fixed estimator $T$, there \emph{always}
exist scalars $a(\theta),b(\theta)$ solving the $L^2$ linear system above (if
the Gram determinant is nonzero), so the $L^2$--projection of $(T-\theta)s$
onto the two--dimensional jet span $T_2(\theta)=\mathrm{span}\{s^{(1)},s^{(2)}\}$ is
well-defined and given by the coefficients in \eqref{eq:ab-sol}.  
\item 
For a Gaussian location model (with variance $\sigma^2$), $\mu_3=0$ and
$\mu_4=3\mu_2^2=3\sigma^4$, and the Gram determinant is nonzero.  The
coefficients simplify; moreover, as we saw above, the first--order
condition already holds for the usual unbiased estimator and it is easily checked that \(b(\theta)=0\), which is why the first--order CRB is tight there.
\item The first--order ODE/projection condition (existence of a scalar
$c_{1,\theta}$ with $(T-\theta)s_\theta=c_{1,\theta}s_\theta^{(1)}$) is
generically not satisfied in canonical exponential families (except in
special cases such as the Gaussian location family where, for instance, \(\mu_3\) vanishes).  This is the statistical
manifestation of a nonzero first--order torsion/vertical component.
\end{enumerate}

\paragraph{Summary}
Our examples thus show that: (i) the Gram entries and hence the canonical projections are expressible in terms of moments of $P_\theta$; (ii) when higher central moments (skewness, kurtosis, etc.) are nonzero, the vertical component of the jet velocity is generically nonzero and contributes to curvature-aware corrections; (iii) only in special cases (notably location families with symmetric kernels, such as Gaussian location) do the torsion components vanish at low order, leading to simple first--order ODE characterisations of efficiency.

\section{Conclusion and Future Work}
\label{sec:future}
This paper developed a Cartan–geometric foundation for the curvature-aware
variance bounds introduced in our earlier Hilbert space formulation of the
Cramér–Rao and Bhattacharyya-type inequalities.  
By representing the square root map 
$s_\theta=\sqrt{f(\cdot;\theta)}$ as a section of the statistical bundle
$E=\Theta\times L^{2}(\mu)$ and analysing its prolongations in the jet tower
$J^{m}(E)$, we showed that efficiency of unbiased estimators in the hierarchy of variance bounds admits a natural interpretation in terms of the hierarchy of holonomy and integrability
conditions associated with Cartan’s contact structure.  

At the analytic level, $m$-th order efficiency is characterised by the
membership condition $(T-\theta)s_\theta \in T_m(\theta)$, where
$T_m(\theta)=\operatorname{span}\{s_\theta^{(1)},\dots,s_\theta^{(m)}\}$.  
We established that this condition is equivalent to the existence of an
$m$-th order linear differential relation satisfied by $s_\theta$, and hence to
the requirement that the prolonged section $j^m s$ take values in the
corresponding ODE-defined submanifold of $J^m(E)$.  
The obstruction to integrability---encoded by the vertical component of the
canonical Ehresmann connection $\nabla^{m}$ on
$\pi_{m,m-1}:J^{m}(E)\to J^{m-1}(E)$---provides a geometric representation of
the inefficiency term appearing in the variance decomposition.  
A model-dependent operator  
$A_\theta: V_{j^ms(\theta)}(J^m/J^{m-1}) \to T_m(\theta)^\perp$  
links this vertical jet component to the Hilbert space residual that enters the
curvature-corrected bound on the variance.  
Thus the Cartan framework unifies the analytic projection picture, the
extrinsic curvature picture in $L^2(\mu)$, and the differential equation viewpoint
on statistical efficiency.

Several directions naturally follow from this work.

\begin{itemize}
\item \textbf{Multi-parameter extensions.}  
  For vector-valued parameter spaces, the torsion and curvature tensors of the
  canonical jet connection interact nontrivially with the Levi–Civita
  connection of the Fisher–Rao metric.  
  This raises the possibility of tensorial generalisations of curvature-aware
  variance bounds and higher-order efficiency conditions.

\item \textbf{Holonomy moduli and geometric control aspects.}  
  The characterisation of estimator efficiency as a holonomy/integrability
  condition suggests studying the moduli space of holonomic versus
  non-holonomic statistical sections.  
  This may link statistical efficiency to geometric control theory and the
  structure of admissible lifts in the jet tower.

\item \textbf{Numerical and computational directions.}  
  In models where $s_\theta^{(m+1)}$ admits closed-form
  expressions (e.g.\ exponential families or Gaussian models),  
  one could possibly derive practical higher-order curvature-corrected bounds in
  high-dimensional or machine learning settings.

\item \textbf{Relation to intrinsic information geometry.}  
  The appearance of jet torsion and non-integrability alongside the Hilbert
  space extrinsic picture suggests a common geometric origin for curvature
  effects arising in both the Fisher–Rao Levi–Civita geometry and the
  $\alpha$-connections of Amari’s dualistic framework.  
  A unified treatment of intrinsic and extrinsic curvature corrections within a
  Cartan–jet setting appears within reach.

\end{itemize}

These directions point toward a broader differential-geometric theory of
statistical efficiency grounded in the geometry of jet bundles, their canonical
connections, and the analytic structure of statistical models.
%%%%%%%%%%%%%%%%%%%%%%%%%%%%%%%%%%%%%%%%%%%%%%%%%%%%%%%%%%%%%

%% if your bibliography is in bibtex format, uncomment commands:
    % Bibliography file (usually '*.%% References with BibTeX database:

\bibliographystyle{elsarticle-num}
\bibliography{ref}






% \end{thebibliography}

\end{document}